\documentclass[12pt]{article}

\usepackage{amsmath}
\usepackage{amssymb}
\usepackage{float}
\usepackage{setspace}
\usepackage[english]{babel}
\usepackage[utf8]{inputenc}
\usepackage{amsthm}
\usepackage[makeroom]{cancel}
\usepackage{wasysym}
\usepackage{verbatim}
\usepackage{graphicx}
\usepackage[margin = 1in]{geometry}
\usepackage[colorinlistoftodos]{todonotes}
\theoremstyle{definition}
	
\theoremstyle{plain}
	\newtheorem{theorem}{Theorem}
    \newtheorem{lemma}[theorem]{Lemma}

\begin{document}
\author{Shilin Ma (Carleton College)\\ Kevin J.~McGown (California State University, Chico)\\ Devon Rhodes (California State University, Chico)\\ Mathias Wanner (Villanova University)}
\title{On the number of primes for which a polynomial is Eisenstein}
\maketitle
\begin{center}
\end{center}
\vspace{2ex}

\section{Introduction}
For an integer $d\geq 2$,
let $f(x) = a_d x^d + a_{d-1} x^{d-1} +\cdots +a_1 x + a_0$ be a polynomial with integer coefficients.
We say that $f$ is Eisenstein if there exists a prime $p$ such that
$p\mid a_i$ for $i = 0,1,\ldots, d-1$, $p^2\nmid a_0$, and $p\nmid a_d$.
The well-known fact that Eisenstein polynomials are irreducible 
is often encountered in an undergraduate algebra course.
See~\cite{bib:C} for a fascinating history of this result, which was proved independently by Sch\"onemann and Eisenstein.

Dobbs and Johnson (see~\cite{bib:DJ}) posed some probabilistic questions concerning Eisenstein polynomials.
In particular, one could ask:
What is the probability that a randomly chosen polynomial is Eisenstein?
Dubickas answers this question in~\cite{bib:D} by providing an asymptotic expression for the number of monic Eisenstein polynomials
of fixed degree and bounded height.  Later Heyman and Shparlinski (see~\cite{bib:HS}) gave an asymptotic expression for the number of Eisenstein polynomials (monic or not) of fixed degree and bounded height but with a stronger error term.  We mention in passing that there are generalizations and variations one may consider; some results in this area include~\cite{bib:H,bib:HS2,bib:MS,bib:DM}.

Our paper builds naturally on~\cite{bib:HS} so we begin by stating their result.
Define the height of a polynomial $f$ to be $\max\{|a_0|,|a_1|,\ldots,|a_d|\}$.
Let $\mathcal{F}_d(H)$ be the set of Eisenstein polynomials of degree $d$ and height at most $H$.
\begin{theorem}[Heyman--Shparlinski]\label{T:HS}
We have
\[
\#\mathcal{F}_d(H)
=\gamma_d(2H)^{d+1}+\begin{cases}
O(H^{d})&\text{if $d>2$}\\
O(H^2(\log H)^2)&\text{if $d=2$}\,.
\end{cases}
\]
\end{theorem}

Let $\psi(f)$ denote the number of primes for which $f$ is Eisenstein.
Our aim is to study the statistics of this function.
We establish the following result, which gives an expression for the mean and variance of the function $\psi(f)$ as $f$ ranges over all
Eisenstein polynomials of a fixed degree.
\begin{theorem}\label{T:1}
If
\[
\alpha_d:=\sum_{p\textrm{ prime}}\frac{(p-1)^2}{p^{d+2}}
\,,\qquad
\beta_d:=\sum_{p\textrm{ prime}}\left(\frac{(p-1)^2}{p^{d+2}}\right)^2
\,,\qquad
\gamma_d:=1-\prod_{p\textrm{ prime}}\left(1-\frac{(p-1)^2}{p^{d+2}}\right)
\,,
\]
then we have
\begin{align*}
&\mu_d:=
\lim_{H\to\infty} \frac{\sum_{f\in \mathcal{F}_d(H)}\psi (f)}{\sum_{f\in \mathcal{F}_d(H)}1}
=
\frac{\alpha_d}{\gamma_d}
\,,
\\[2ex]
&\sigma_d^2:=
\lim_{H\to\infty} \frac{\sum_{f\in \mathcal{F}_d(H)}(\psi(f)-\mu_d)^2}{\sum_{f\in \mathcal{F}_d(H)}1}
=
\frac{\alpha_d+\alpha_d^2-\beta_d-\mu_d\alpha_d}{\gamma_d}
\,.
\end{align*}
\end{theorem}
We note in passing that $\alpha_d$ and $\beta_d$ can be expressed as finite linear
combinations of values of the prime zeta function $P(s)=\sum_{p} p^{-s}$.
Throughout this paper, the variables $p$ and $q$ will always denote primes.
See Section~\ref{S:constants} for additional comments on $\alpha_d$, $\beta_d$, $\gamma_d$, $\mu_d$, $\sigma_d^2$,
including a table of numerical values for various values of $d$.

\section{Proofs}
As usual we let $\omega(n)$ denote the number of distinct prime factors of $n$
and let $\phi(n)$ denote the Euler phi-function.
Following~\cite{bib:HS}, we let $\mathcal{H}_d(s,H)$ be the number of polynomials of degree $d$ and height at most $H$
satisfying
$s\mid a_i$ for $i = 0,1,\ldots, d-1$,
$\gcd(a_0/s,s)=1$, and
$\gcd(a_d,s)=1$.


\begin{lemma}\label{L:HS}
We have
\begin{equation}
\#\mathcal{H}_d(s,H) = \frac{(2H)^{d+1}\phi^2(s)}{s^{d+2}}+O\left(\frac{2^{\omega(s)}H^d}{s^{d-1}}\right)
\,.
\end{equation}
\end{lemma}

\begin{proof}
See Lemma~5 of~\cite{bib:HS}.
\end{proof}

\begin{lemma}\label{L:1}
We have
\begin{equation}\label{E:1}
\sum_{f\in\mathcal{F}_d(H)}\psi (f)=(2H)^{d+1}\alpha_d + \begin{cases}
O(H^2) & \text{if $d>2$}\\
O(H^2\log\log H) & \text{if $d=2$}\,.
\end{cases}
\end{equation}
\end{lemma}

\begin{proof}
We rewrite the sum in question as a sum over primes and apply Lemma~\ref{L:HS};
this yields
\begin{align*} 
\sum_{f\in \mathcal{F}_d(H)}\psi (f) &=
\sum_{p\leq H} \#\mathcal{H}_d(p,H)\\ &= \sum_{p\leq H}
\left[\frac{(2H)^{d+1}\phi^2(p)}{p^{d+2}}+O\left(\frac{2^{\omega(p)}H^d}{p^{d-1}}\right)\right] \\
&=(2H)^{d+1}\sum_{p\leq H} \frac{(p-1)^2}{p^{d+2}}+\sum_{p\leq H} O\left(\frac{H^d}{p^{d-1}}\right)\\
&=(2H)^{d+1}\sum_{p} \frac{(p-1)^2}{p^{d+2}}-(2H)^{d+1}\sum_{p>H} \frac{(p-1)^2}{p^{d+2}}+\sum_{p\leq H} O\left(\frac{H^d}{p^{d-1}}\right)
\,.
\end{align*}
The splitting of $\sum_{p\leq H}$ into $\sum_p$ and $\sum_{p>H}$ is justified since $\sum_p$
converges absolutely.
It remains to bound the second and third terms in the last line above.
We bound the second term using the integral test to obtain
\begin{equation*}
(2H)^{d+1}\sum_{p>H} \frac{(p-1)^2}{p^{d+2}} = O\left(H^{d+1}\int_{H}^{\infty}\frac{(x-1)^2}{x^{d+2}}\;dx\right)=O\left(H^{d+1}H^{-d+1}\right)=O\left(H^2\right)
\,.
\end{equation*}
For the third term, we find
\[
H^d \sum_{p\leq H} \frac{1}{p^{d-1}}
=
\begin{cases}
O\left(H^2\right)&\text{if $d>2$}\\
O\left(H^2\log\log H\right)&\text{if $d=2$}\,,
\end{cases}
\]
where we have used Mertens' Theorem (see, for example,~\cite{bib:P}) in the case of $d=2$.
\end{proof}

%

\begin{lemma}\label{L:2}
We have
\begin{equation*}
\sum_{f\in\mathcal{F}_d(H)}\psi(f)^2=
(\alpha_d+\alpha_d^2-\beta_d)(2H)^{d+1}
 + \begin{cases}
O(H^2) & \text{if $d>2$}\\
O(H^2(\log\log H)^2) & \text{if $d=2$}\,.
\end{cases}
\end{equation*}
\end{lemma}

\begin{proof}
If we define
\begin{align*}
\tau(f,p) = \begin{cases}
1& \text{if $f$ is $p$-Eisenstein} \\
0 & \text{otherwise}\,,
\end{cases}
\end{align*}
then the first sum can be rewritten as
\begin{align*}
\sum_{f\in\mathcal{F}_d(H)} \psi(f)^2 &= \sum_{f\in\mathcal{F}_d(H)}\left(\sum_{p\text{ prime}} \tau(f,p)\right)^2 \\
&=\sum_{f\in\mathcal{F}_d(H)}\left(\sum_{p \text{ prime}}\tau(f,p)\sum_{q\text{ prime}}\tau(f,q)\right) \\
&=\sum_{f\in\mathcal{F}_d(H)}\left(\sum_{p,q\text{ prime}}\tau(f,p)\tau(f,q)\right) \\
&=\sum_{p,q\text{ prime}}\left(\sum_{f\in\mathcal{F}_d(H)}\tau(f,p)\tau(f,q)\right) 
\,.
\end{align*}

The inner sum above represents the number of polynomials of height at most $H$ that are Eisenstein for both $p$ and $q$, but the fact that $p$ may equal $q$ complicates matters. Consequently, we have
\begin{align*}
\sum_{f\in\mathcal{F}_d(H)} \psi(f)^2 &=
\sum_{\substack{p\leq H}}
\#\mathcal{H}(p,H)
+
\sum_{\substack{pq\leq H\\p\neq q}}
\#\mathcal{H}(pq,H)
\,.
\end{align*}
The first sum on the right-hand side above is exactly what appears in Lemma~\ref{L:1},
and therefore is it equal to the right-hand side of (\ref{E:1}).
It remains to deal with the second sum, which equals
\begin{align*}
&
\sum_{\substack{pq \leq H\\p\neq q}} \#\mathcal{H}(pq,H)\\
&=
(2H)^{d+1}\sum_{\substack{pq \leq H\\p\neq q}} \frac{(p-1)^2(q-1)^2}{p^{d+2}q^{d+2}} + O\left(\sum_{\substack{p,q\text{ prime}\\pq \leq H}}\frac{H^d}{(pq)^{d-1}}2^{\omega(pq)}\right)
\\
&=
(2H)^{d+1}\sum_{\substack{pq \leq H}} \frac{(p-1)^2(q-1)^2}{p^{d+2}q^{d+2}}
-
(2H)^{d+1}\sum_{\substack{p^2 \leq H}} \left(\frac{(p-1)^2}{p^{d+2}}\right)^2
+ O\left(H^d\sum_{\substack{p,q\text{ prime}\\pq \leq H}}\frac{1}{(pq)^{d-1}}\right)
\,.
\end{align*}
For the first term, as in the proof of Lemma~\ref{L:1}, we have
\begin{align*}
(2H)^{d+1}\sum_{\substack{pq \leq H}} \frac{(p-1)^2(q-1)^2}{p^{d+2}q^{d+2}}
&=
(2H)^{d+1}\sum_{p,q} \frac{(p-1)^2(q-1)^2}{p^{d+2}q^{d+2}}
+
(2H)^{d+1}\sum_{p>H}\frac{1}{p^d}\sum_{q>H/p}\frac{1}{q^d}
\\
&=
(2H)^{d+1}\left(\sum_{p} \frac{(p-1)^2}{p^{d+2}}\right)^2
+O\left(H^{d+1}\sum_{p>H}\frac{1}{p^d}\right)
\\
&=
(2H)^{d+1}\alpha_d^2
+O(H^2)
\,.
\end{align*}
For the second term,
\begin{align*}
(2H)^{d+1}\sum_{\substack{p^2 \leq H}} \left(\frac{(p-1)^2}{p^{d+2}}\right)^2
&=
(2H)^{d+1}\sum_{p} \left(\frac{(p-1)^2}{p^{d+2}}\right)^2
-
(2H)^{d+1}\sum_{\substack{p> \sqrt{H}}} \left(\frac{(p-1)^2}{p^{d+2}}\right)^2
\\
&=
(2H)^{d+1}\beta_d+O(H^{3/2})
\,.
\end{align*}
Finally, for the third term, we have
\begin{align*}
H^d
\sum_{\substack{p,q\text{ prime}\\pq \leq H}}\frac{1}{(pq)^{d-1}}
=
\begin{cases}
O(H^2), & \text{if $d>2$}\\
O(H^2\left(\log\log H \right)^2), & \text{if $d=2$}\,.
\end{cases}
\end{align*}
Putting this all together proves the lemma.
\end{proof}


\begin{proof}[Proof of Theorem~\ref{T:1}]
The part of the theorem concerning the mean $\mu_d$ follows immediately from Lemma~\ref{L:1} and Theorem~\ref{T:HS}.  Now we consider the variance:
\begin{align*}
\sigma^2_d
&=
\lim_{H\to\infty} \frac{\sum_{f\in \mathcal{F}_d(H)}(\psi(f)-\mu_d)^2}{\sum_{f\in \mathcal{F}_d(H)}1}\\
&=
\lim_{H\to\infty}
\frac{1}{\#\mathcal{F}_d(H)}
\sum_{f\in\mathcal{F}_d(H)}\left(\psi(f)^2 -2 \psi(f)\mu_d +\mu_d^2\right)\\
&=
\lim_{H\to\infty}\frac{1}{\#\mathcal{F}_d(H)}\left[\sum_{f\in\mathcal{F}_d(H)}\psi(f)^2-2\mu_d\sum_{f\in\mathcal{F}_d(H)}\psi(f)+\mu_d^2\sum_{f\in\mathcal{F}_d(H)}1\right]
\,.
\end{align*}
By Lemma~\ref{L:1}, Lemma~\ref{L:2}, and Theorem~\ref{T:HS}, the limit above equals
\[
  \frac{1}{\gamma_d}\left[
  (\alpha_d+\alpha_d^2-\beta_d)
  -2\mu_d\alpha_d
  +\mu_d^2\gamma_d
  \right]
  \,,  
\]
which simplifies to the desired expression.
\end{proof}

\section{Remarks on the constants}\label{S:constants}

It is not hard to show that
\[
  \alpha_d=\frac{1}{2^{d+2}}+O\left(\frac{1}{3^d}\right)
  \,,\qquad
  \beta_d=\frac{1}{2^{2(d+2)}}+O\left(\frac{1}{3^{2d}}\right)
  \,,\qquad
  \gamma_d=\frac{1}{2^{d+2}}+O\left(\frac{1}{3^d}\right)
  \,.
\]
It then follows that $\lim_{d\to\infty}\mu_d= 1$ and $\lim_{d\to\infty}\sigma_d^2= 0$, as one would expect.
If one was interested in the mean $\hat{\mu_d}$ and variance $\hat{\sigma}_d^2$ of $\psi(f)$ as $f$ ranges over all polynomials, instead of just Eisenstein polynomials, one would obtain the simpler expressions
$\hat{\mu}_d=\alpha_d$ and $\hat{\sigma}_d^2=\alpha_d-\beta_d$.
We will not prove this explicitly but it essentially follows from the proof of Theorem~\ref{T:1}.
In this case, one observes that $\lim_{d\to\infty}\hat{\mu_d}=0$ and $\lim_{d\to\infty}\hat{\sigma}_d^2=0$,
as expected.


\begin{table}[h]
\centering
\begin{tabular}{c|c|c|c|c|c|c}
$d$ & $\alpha_d=\hat{\mu}_d$ & $\beta_d$ & $\gamma_d$ & $\mu_d$ & $\sigma^2_d$ & $\hat{\sigma}_d^2$ \\
\hline
$2$&$0.17971$&$0.00731$&$0.16765$&$1.07192$&$0.07187$&$0.17239$\\ 
$3$&$0.05653$&$0.00127$&$0.05557$&$1.01714$&$0.01705$&$0.05525$\\
$4$&$0.02255$&$0.00027$&$0.02243$&$1.00519$&$0.00517$&$0.02227$\\
$5$&$0.00989$&$0.00006$&$0.00988$&$1.00169$&$0.00169$&$0.00983$\\
$6$&$0.00456$&$0.00001$&$0.00456$&$1.00056$&$0.00056$&$0.00454$\\

\end{tabular}
\caption{Approximate values of the constants for small $d$}
\end{table}

\section*{Acknowledgement}
This research was completed as part of the Research Experience for Undergraduates and Teachers program at California State University, Chico funded
by the National Science Foundation (DMS-1559788).

\end{document}